\documentclass[11pt]{article}

\usepackage{amsmath, epsfig, cite}

\usepackage{amssymb}

\usepackage{amsfonts}

\usepackage{latexsym}

\newcommand{\old}[1]{{}}

\newtheorem{theorem}{Theorem}[section]

 \numberwithin{equation}{section}

 \makeatletter \@addtoreset{equation}{section} \makeatother

\newcommand\blfootnote[1]{%
  \begingroup
  \renewcommand\thefootnote{}\footnotetext{#1}%
  \addtocounter{footnote}{-1}%
  \endgroup
}

\title{Minimizing the number of episodes and Gallai's theorem on intervals}

\author{\'Eva Czabarka$^{2,3,5}$,
\\ University of South Carolina, 
{\tt czabarka@math.sc.edu}\\
L\'aszl\'o  A. \ Sz\'ekely$^{1,2,3,4,5}$
\\
University of South Carolina, {\tt szekely@math.sc.edu}\\
Todd Vision$^{1,3}$
\\
University of North Carolina, {\tt tjv@biol.unc.edu}\\
}
\blfootnote{\\
$^1$ This work started at the Phylogeny program of the Isaac Newton Institute, Cambridge.\\
$^2$ This author was supported in
part by  a Marie
Curie Fellowship HUBI MTKD-CT-2006-042794.\\
$^3$ This author was supported in
part by the NIH NIGMS contract  1 R01 GM078991-01 and 3 R01 GM078991-03S1.\\
$^4$ This author was supported in
part by the NSF DMS
contracts 0701111 and 1000475.\\
$^5$ This author was supported in
part by DARPA and AFOSR under the contract FA9550-12-1-0405.
}

\begin{document}

\maketitle

\begin{abstract}
In 1996, Guigo et al. [Mol. Phylogenet. Evol., {\bf 6} (1996), 189--203] posed the following problem: for a given species
tree and a number of gene trees, what is the minimum number of duplication
episodes, where several genes could have undergone duplication together to generate
the observed situation. (Gene order is neglected, but duplication of
genes could have happened only on certain segments that duplicated).
We study two versions of this problem, one of which was algorithmically solved not long ago
by Bansal and Eulenstein \cite{eulenstein}. We provide min-max theorems
for both versions that generalize Gallai's archetypal min-max theorem
on intervals, allowing simplified proofs to the correctness of the algorithms (as it always
happens with duality)  and deeper understanding.  An interesting feature of our approach
is that its recursive nature requires a generality that bioinformaticians attempting to solve
a particular problem usually avoid. 
\end{abstract}

\eject


\section{Introduction}
In 1996 Guigo et al. \cite{guigo} 
posed the following problem: for a given species
tree and a number of gene trees, what is the minimum number of episodes
of gene duplication, where several genes could have duplicated in any
single episode. The constraints of the problem include some vertices
of the gene trees identified as duplication vertices; and duplication
vertices have some associated intervals in the species tree, where a
duplication  of the gene represented by the gene tree
must have taken place. We give mathematical definitions in Section~\ref{combin},
and explain the relevance of our results in Section~\ref{bioinf}. 

Several variants of this problem have been investigated: 
 \cite{fellows}, \cite{vertebrate}, \cite{speedup} , \cite{infer}, \cite{luo}, \cite{ilp}.
Bansal and Eulenstein \cite{eulenstein}
solved a version  of this long-standing open problem with a greedy algorithm
and proved the correctness of the algorithm by induction.

The purpose of our note is to put these problems and  results into  proper
combinatorial context. There is no need to assume that the trees have no internal vertices
of degree two or that they 
are binary. The intervals associated with the duplication 
vertices can be defined differently 
from the definitions  in the biology literature.
The greedy algorithm still
works, furthermore,  simple min-max theorems give a good
characterization to the minimum number of episodes, even in this
more general setting. As usual, the duality allows for more transparent proofs for the 
correctness of the optimization algorithms.

The min-max theorems are straightforward
generalizations to Gallai's Theorem on intervals  (Gallai did not 
publish actually this theorem, and it was first printed 
in a paper of Hajnal and Sur\'anyi \cite{hs}):
 
\begin{theorem} {\rm [Gallai]} Let us be given a finite family of \label{galla}
closed intervals on a straight line. Denote by $\nu$ the size of the
largest set of pairwise disjoint intervals, and by $\tau$ the smallest
number of points that can cover all intervals i.e. every interval contains
at least one of the points. Then $\nu =\tau$ holds.
\end{theorem}

As a reminder, we reproduce here the proof to Gallai's Theorem,
as our proofs were developed to generalize it.
Clearly $\nu \leq \tau$, as disjoint intervals have to be covered
by distinct points. We show that $\nu$ points suffice to cover all
intervals. Apply the following algorithm recursively until the
interval system is empty:

{\it Pick the leftmost right endpoint from  all right endpoints of  
 intervals from the family,
and delete all intervals from the system that this point covered. Add
the picked point to the list of selected points.  }

We picked  some right endpoints of intervals that are pairwise
disjoint by the construction, therefore these endpoints are at most
$\nu $ in number. These right endpoints cover all our intervals by the  
construction.
$\spadesuit$

Gallai's Theorem have been generalized by Sur\'anyi (see \cite{gyarfas}) essentially
with the same proof:
\begin{theorem} {\rm [Sur\'anyi]}   \label{sura} Let us be given a finite family of
subtrees on a tree. Denote by $\nu$ the size of the
largest set of pairwise disjoint subtrees, and by $\tau$ the smallest
number of points that can cover all subtrees i.e. every subtree contains
at least one of the points. Then $\nu =\tau$ holds.
\end{theorem}

\section{Describing the combinatorial problems} \label{combin}
Let us be given a finite set $X$. The elements  of $X$ are called {\em taxa}.
Let us be given a tree $S$ with root $R$, such that the leaves of $S$ 
are labelled with elements from $X$ in a one-to-one manner.
Root $R$ is joined by an edge to $\infty$. We call $S$ the {\em species tree}
for the taxa in $X$. When we talk about vertices of $S$, we exclude  $\infty$.

Let us be given $k$ {\em gene trees}, say  for $i=1,2,...,k$ the gene tree $G_i$.
The leaves of $G_i$ are labelled with some taxa from $X$, but a taxon can occur in  
more than
one leaf, and not all taxa are necessarily represented by a leaf in $G_i$.
A leaf corresponds to one taxon only. 
We assume that $G_i$ also has a root $R_i$ and one more edge
going from the root to $\infty_i$. $\infty_i$ is not considered a vertex of the  
gene tree.

Vertices of $S$ have a natural partial order, namely $u\geq_S v$, if $u=v$ or $u$
separates $v$ from $\infty$.  When we  speak about the interval of $u$ and $v$ in $S$,
we mean the interval in this natural partial order.
We call $u$ the {\em upper endpoint} and $v$ the {\em lower endpoint} 
of this interval. Vertices of  $G_i$ have a natural partial
order $\geq_i$ defined similarly. $>_S$ and $>_i$ will refer to strict inequalities
in these partial orders.

Assume further that every $G_i$ has a subset $D_i$ of its vertices specified
that are called {\it duplication vertices}. For every $i$ and every
$d\in D_i$, we have an associated path $P$,
which is a subpath of a path connecting $\infty$ to a leaf   in the species tree $S$. The ordered pair
$(P,d)$   will be called the
{\it  duplication interval}  associated  to the duplication vertex $d\in D_i$, and for more convenient notation
we write it as $P_d$.
In this way we maintain names on the duplication intervals which tell  which
duplication vertex of which gene tree generated the duplication 
interval. The same  intervals can have multiple names as duplication intervals: the same path  $P$  in $S$
may be assigned as a duplication interval to vertices in different gene trees, and also to
several pairs of $\leq_i$-comparable (or not $\leq_i$-comparable) duplication
vertices of the same $G_i$. Duplication intervals with different names are considered
distinct objects although their underlying  intervals in $S$ are the same.

The following
monotonicity assumption is made on the associated duplication  intervals:
\begin{equation} \label{treelike}
\forall i \forall d,e \in D_i \   \   \  d\geq_i e \rightarrow
(\max_S P_d\geq_S \max_S P_e) \land (\min_S P_d\geq_S \min_S P_e).
\end{equation}

We say that for two duplication intervals  $P_e<_i P_d$,
if $e<_i d$ for the duplication vertices $e,d$ in $G_i$.
A {\it chain} of duplication 
intervals   is a sequence
of duplication intervals associated  to duplication vertices $e_1<_i e_2<_i ...  
<_i e_m$
for some $i=1,2,...,k$. (We freely change between speaking  about chains of duplication intervals
and chains of duplication vertices in the gene trees, as they are in bijective correspondence.)
  We may have several
copies of a path in $S$ present as a duplication interval, and it depends on the  
label
of a particular copy whether it satisfies a $P_e<_i P_d$  type relation
or belongs to a certain chain. Some copy of a path may do it, while
another may not. Also, it may happen that a single path of $S$
satisfies a strict $P_e<_i P_d$ ordering with proper duplication vertices $e,d\in D_i$.

Now we have the following  models and  optimization problems as:\\
{\bf Discrete model.}
Let $V(S)$ denote the vertex set of the species tree $S$. Let $V^*(S)$
denote the extension of $V(S)$ by allowing unlimited number of 
 copies of the vertices. 
We distinguish these copies from each other, but keep the information on
which vertices of  $V^*(S)$ are copies of the same vertex of  $V(S)$.
The elements of $V^*(S)$ inherit
the $\geq_S$ partial order, if they are copies of different vertices,
and are incomparable when they are copies of the same vertex.
We denote this extended partial order by $\geq_{S^*}$.

Consider  maps $f: \cup_{i=1}^k D_i \rightarrow V^*(S)$, which  
have the
property that for all $i$,  the restriction $f_{|D_i}: D_i\rightarrow f(D_i)$ preseves the
partial order  in the following sense:  [$d>_i e$ implies  that $f(d)>_{S^*} f(e)$ or
$f(d)$ and $f(e)$ are different copies of the same vertex]. 
A value
of $f$ is called a {\em (duplication)  episode.} 

{\bf Objective:} Minimize the
quantity $|f( \cup_{i=1}^k D_i)|$ over all 
 maps $f$, and/or characterize optimal solutions.
(A fast algorithm for this minimization 
was discovered  by Bansal and Eulenstein  \cite{eulenstein}.)

We consider an alternative model as well:

{\bf Continuous model.} 
Consider the edges of $S$ as {\em line segments} that have
interior points. Let $int(S)$ denote the union of the set of interior  
points
of all edges of $S$. The partial order $\geq_S$ on $S$ naturally extends to
$\tilde S=V(S)\cup int(S)$. The extension will be denoted by $ {\geq}_{\tilde S} $.

Consider now maps $f: \cup_{i=1}^k D_i \rightarrow \tilde S$, which  
have the
property that for all $i$,  the restriction $f_{|D_i}: D_i\rightarrow f(D_i)$ strictly preserves
 the
partial order  [i.e. 
$d>_i e$ implies  $f(d)>_{\tilde S} f(e)$]. We still call the values of $f$ (duplication)
 episodes.

{\bf Objective}: Minimize the
quantity $|f( \cup_{i=1}^k D_i)|$ over all 
 maps $f$, and/or characterize optimal solutions.

We say that a duplication interval in the species tree is {\em degenerate},
if it has only one point. For the continuous model we require, in  addition,
that {\em the duplication intervals are non-degenerate}, as otherwise the problem
may not have a feasible solution at all. 
Note that the minimum number of episodes in the discrete and continuous
models can be different, even if all duplication intervals are non-degenerate.

\section{The new min-max theorems}
First we discuss the simpler {\em discrete model}. Let $\Lambda$ be an arbitrary
index set. 
For $ \lambda\in\Lambda$, let $C_\lambda $ be a set of duplication 
intervals that make a chain with respect to one of the gene orders. We call $\{C_\lambda: \ 
\lambda \in \Lambda\}$ a {\it disjoint chain packing}, if 
for every
$\lambda\not= \lambda'$, elements of  $C_\lambda $ and $C_{\lambda'} $
 do not share
vertices in $S$, i.e. $(\cup C_\lambda)\cap (\cup C_\lambda')=\emptyset$.
We call $\sum_{\lambda\in\Lambda} |C_\lambda|$ the {\em value} of the disjoint 
chain packing.
Fix now an arbitrary disjoint chain packing,  $\{C_\lambda: \ 
\lambda \in \Lambda\}$.
 Now the  number of episodes needed is clearly at least
as much as the value of this chain packing, as different members
of a chain must belong to different episodes and 
vertex disjoint chains must use disjoint sets of episodes.
\begin{theorem} \label{Disc}
In the discrete model, the minimum number of duplication episodes
equals to the maximum value of   a disjoint chain packing.
\end{theorem} 
We will prove the other (non-trivial) inequality in the next section.

We continue with the  {\em continuous model}.
For every $\lambda\in\Lambda$, let $C_\lambda $ be a set of duplication 
intervals that make a chain with respect to one of the gene orders. We call 
$\{C_\lambda: \ \lambda \in \Lambda\}$ 
an {\it almost disjoint chain packing}, if 
the following restrictions for
 intersections  (in $S$) of elements from different chains,   $C_\lambda $ and $C_\lambda' $, hold:\\
 (i) for any $U\in C_\lambda $ and $U'\in C_\lambda' $, we have $|U \cap U'|\leq 1$;\\
 (ii) for any $U\in C_\lambda$ and $U'\in C_{\lambda'}$,
$|U \cap U'|= 1$ imply that the single element of $ U \cap U'$ is the 
 $>_S$ (upper)  endpoint of at least one of $U$ and $U'$;\\ 
(iii) if 
$U\in C_\lambda$ and $U'\in C_{\lambda'}$ intersect in 
a single point that is the $>_S$ endpoint  of $U$,  but not the  $>_S$ endpoint  of  $U'$,
then there is an $R\in C_{\lambda'}$, such that $U'>_{\lambda'} R$,    $U\cap R= U\cap U'$, and
this common intersection point is the    $>_S$ (upper)  endpoint of $R$ as well.\\ 
Of course,  duplication intervals from the same chain are allowed to intersect.

Note that condition (iii) means that different
chains from an almost disjoint chain packing, as sets in $\tilde S$, may only intersect
at {\em nodes} of $S$, and if $v$ is a node where several chains intersect,
then it is the $>_S$ upper endpoint of all the chains that go through it with
{\em at most one exception}. The chains that go through $v$
all go down along different edges from $v$, and the exceptional
chain must contain an interval that has $v$ as its $>_S$ upper endpoint.

 For $v\in V(S)$, and an almost disjoint chain packing 
$\{C_\lambda: \ \lambda \in \Lambda\}$,
let ${\cal E}_\Lambda(v)$ denote the number of chains $C_\lambda$,
which have elements with upper
endpoint $v$. Fix  an arbitrary
almost disjoint chain packing, $\{C_\lambda: \ \lambda \in \Lambda\}$.
We call 
 $$\sum_{\lambda\in\Lambda} |C_\lambda|-\sum_{v\in V(S), {{\cal E}(v)\geq 1}}
( {\cal E}_\Lambda(v)-1) $$ the {\em value} of the almost disjoint chain packing.
Now the  number of duplication episodes needed is clearly at least
the value of the almost disjoint chain packing, as  different members
of any chain must belong to different episodes, disjoint intervals also must
belong to different episodes, 
and for any vertex $v$ with ${\cal E}_{\lambda}\ge 1$, we may use $v$ as 
the episode for (no more than) one of the intervals
from each of the ${\cal E}_{\lambda}(v)$
chains covering $v$. Thus,
we can save on the duplication intervals containing  $v$, by using $v$,
${\cal E}_\Lambda(v)-1$ episodes, compared to not using $v$ as an episode. 
\begin{theorem} \label{Cont}
In the continuous model, the minimum number of duplication episodes
equals to the maximum  value of almost disjoint chain packings.
\end{theorem} 
We will prove the other (non-trivial) inequality in the next section.

It is easy to see that  Gallai's Theorem~\ref{galla} is a special instance
of both Theorems \ref{Disc} and \ref{Cont}, when the species tree
is a path (only one taxon is present) and every gene tree
has a single duplication event.

\section{Proofs}
As the algorithms and the proof of their correctness through the
respective min-max theorems are very similar, we describe them 
in one text, and tell the differences as they arise.

The proof is mathematical induction on the total number of 
duplication vertices in the gene trees.
There is nothing to prove if none of the gene trees contain any
duplication vertex, and in this case the empty (almost) disjoint chain 
packing suffices.
The algorithm will remove the duplication
designation of certain vertices in the gene trees, but not the vertices
themselves; and will solve recursively the reduced  problem with the reduced
number of duplication vertices.
We will also provide (almost) disjoint chain packing for the reduced problem,
with the right value,
 such that the min-max theorem holds
by the inductive hypothesis for the reduced problem. Then, case by case,
we show that the number of episodes from the reduced problem plus
the number of episodes created by our greedy algorithm  in the reduction step
equals to the size of an (almost) disjoint chain packing for the original
problem. This will show simultaneously the optimality of our 
greedy algorithm and the truth of the corresponding min-max theorem.  

So we assume that we already know that the recursive algorithm solves the problem
optimally in any instance when 
 the total number of duplication vertices is less than the current amount and
 that in these instances a  disjoint/almost disjoint chain packing can also be built with value equal to the
  minimum number of episodes..

Every duplication interval has an $\leq_S$-upper endpoint.
Find a $\leq_S$-minimal among  all $\leq_S$-upper duplication 
interval endpoints. Let this vertex of $S$ be $P$. 

\noindent{\bf Discrete model}: 
Let $k\geq 1$ be the largest integer such that  $P$ is
$\leq_S$-upper endpoint of each of the $k$ elements of some chain 
$<_j$ for the  order in a gene tree $G_j$, say
 $L_1<_jL_2<_j\cdots <_j L_k$. Remove  
the duplication designation
of any vertex $d$ in any  gene tree $G_i$, if  $P$ belongs to the
duplication interval of $d$ and no $<_i$-chain of duplication
vertices in $G_i$ with maximum element $d$ has length $ k+1$.  

By induction, the same recursive
algorithm solves the reduced episode problem optimally such that
the min-max theorem holds for the reduced problem. 
Add $P$ with multiplicity
$k$  to the system of episodes. We  construct  recursively a disjoint chain packing providing 
the same value,
through the following two cases:

(i) no chain in the  optimal disjoint chain packing for the  
 reduced problem covers $P$.  Add  the chain $\{L_1<_jL_2<_j \cdots <_jL_k\}$ to the 
 disjoint chain packing for the reduced problem---note that we still have a disjoint chain packing.
We have the following chain of inequalites:
 the minimum number of episodes in the original problem is at most the 
minimum number of episodes
in the reduced problem $+k$, which equals to the maximum value of a  disjoint chain packing in the
reduced problem $+k$, which is at most the maximum value of a disjoint 
chain packing in the original 
problem. We already know the trivial inequality for the min-max
theorem, hence in this case our algorithm provides  the same number of episodes
as the value of a disjoint chain packing.

(ii) a chain ${C}$ in the optimal   disjoint chain packing
 for the
 reduced problem covers $P$. Let the lowest element in the  chain ${C}$ 
correspond to the duplication interval $U$.  By the choice of $P$, $P$ must be in $U$.
The duplication vertex $d$, which is responsible for $U$,
has not been deleted   from the list of duplication vertices. 
This means that $d$ is the maximum element in a
$(k+1)$-chain  ${ C}'$ of duplication vertices in his gene tree. 
Merge ${C}$ and ${C}'$ into a single chain (it is possible as $d$ was lowest element $C$
but highest in $C'$), and add
the merged chain to the optimal disjoint chain packing for the reduced problem to obtain
a disjoint chain packing for the original problem. The number of episodes that we use
for the original problem equals to the value of the disjoint chain packing that we constructed
for the original problem.

In both cases, we constructed a disjoint chain packing, whose value is the same
as the number of episodes constructed, and hence the induction proof is complete.

\noindent{\bf Continuous model}: 
Note that $P$ is not a leaf vertex in $S$, as duplication intervals
in the continuous model are non-degenerate.
Assume that  $e_1$, $e_2$, ..., $e_\ell$ are the edges leaving $P$
in directions different from $\infty$ in $S$. 
For $j=1,2,...,\ell$,
let ${H}_j $ denote a
longest chain of duplication intervals over all  gene trees with
the following properties:\\ 
($\alpha$) the upper endpoint of every duplication  interval
from the chain is $P$.\\ 
($\beta$) every duplication interval of the chain uses the edge $e_j$. \\
Set $d_j= |H_j|$. As we will not need the $d_j=0$ terms,
assume that only those edges leaving $P$ are enumerated on which
$d_j\geq 1$, and for convenience those edges are still labelled as 
$1,2,...,\ell$.

We create a reduced problem by removing the duplication vertex designation
of certain vertices in the gene trees. Assume that $I_f$ is a duplication
interval containing $P$, with vertex $f\in D_i$ from the gene tree $G_i$.
  We remove  the duplication designation of $f$
and the duplication interval $I_f$ if\\
(a) $P$ is the upper endpoint of $I_f$, or\\
(b) $P$ is a vertex of $I_f$ but not an endvertex (so $I_f$ uses 
some $e_j$ edge from $P$) 
and in $G_i$,  there is no chain of duplication vertices of length exceeding
$d_j$ in which $I_f$ is the top element, or\\
(c)  $P$ is the lower endpoint of $I_f$, and for every $j=1,2,...,\ell$,
 in $G_i$ there is no chain of duplication vertices of length exceeding
$d_j$ in which $I_f$ is the top element and all other elements have 
duplication intervals passing through $e_j$.

By mathematical induction, the same recursive algorithm solves the reduced episode problem 
optimally such that the min-max theorem holds for the reduced problem.
Add to the list of episode locations
 the following points  in $\tilde S$:
$P$ itself and $d_j-1$ distinct points from the interior of $e_j$ for  every $j=1,2,...,\ell$.
We  construct recursively an almost disjoint 
chain packing providing the same value, through the following two cases:

(i) no chain in the 
  optimal almost disjoint chain packing 
 for the reduced problem covers $P$. Add to this 
almost disjoint system of chains, which provides the min-max result for the
reduced episode problem by hypothesis,
  a length $d_j$ chain  for every $j=1,2,...,\ell$ from a gene tree, such that
every duplication interval of this  length $d_j$ chain has upper endpoint $P$ and 
uses the edge $e_j$. It is easy to see that we obtained an almost
disjoint chain packing for the original problem. Simple calculation
shows, analogously to the discrete case,  that from the min-max result
for the reduced problem, we obtain that in original problem  the number of episodes
equals to the value of the almost disjoint chain packing,
as both sides increase by $1+\sum_{j=1}^\ell (d_j-1)=\big(\sum_{j=1}^{\ell} d_j\big)-(\ell-1)$. 
 
The alternative of (i) is that one or more chains in the 
  optimal almost disjoint chain packing 
 for the reduced problem covers $P$.  Observe that $P$ may belong to two chains only if $P$ is upper endpoint
 of some elements of one of the  chains. However, (a) has removed those elements. We are left with: 

(ii) a single chain $C$ (corresponding to some gene tree $G_i$)
in the optimal almost disjoint chain packing for the
 reduced problem covers $P$. 
Let the smallest element in this chain 
be $U$,  so $P\in U\in C$ by the choice of $P$. 

The duplication vertex $d$ that is responsible for $U=U_d$
has not been deleted.
According to the removal rules, $P$ cannot be 
the upper endpoint of $U$. Therefore $P$ is an internal point or lower endpoint of $U$.\\
 If $P$ is an internal point of $U$  
 and $U$ passes through $e_j$, then there is a 
chain $C'$ of duplication vertices in the gene tree corresponding to $U$,
in which $U$ is the $(d_j+1)^{th}$ element. 
Create the almost disjoint system of chains for the original problem
from the optimal almost disjoint system of chains for the reduced problem in the
following way:
replace $C$ with $C\cup C'$, which is still a chain; and for 
$t=1,2,...,\ell, t\not= j$, add the chain ${H}_t$. Indeed, we
obtain an almost disjoint chain packing. Simple calculation
shows, analogously to the discrete case,  that from the min-max result
for the reduced problem, we obtain 
that the number of episodes equals to the value of the almost disjoint chain packing
in the original problem, as 
both sides increase by $1+\sum_{j=1}^\ell (d_j-1)=\big(\sum_{j=1}^{\ell} d_j\big)-(\ell-1)$. \\
 If $P$ is the lower endvertex of the duplication interval $U$  that comes from the gene tree
 $G_i$, then for some $j$, there is
a $(d_j+1)$-chain   $C'$ in $G_i$ with top element $U$, and all other
elements of this chain use $e_j$ in their duplication intervals. 
Create the almost disjoint system of chains for the original problem
from the optimal almost disjoint system of chains for the reduced problem in the
following way:
replace $C$ with $C\cup C'$, which is still a chain; and for 
$t=1,2,...,\ell, t\not= j$, add the chain ${H}_t$. Indeed, we
obtain an almost disjoint chain packing. Simple calculation
shows, analogously to the discrete case,  that from the min-max result
for the reduced problem, we obtain 
that the number of episodes equals to the value of the almost disjoint chain packing
in the original problem, as 
both sides increase by $1+\sum_{j=1}^\ell (d_j-1)$. \\

Cases (i) and (ii) together will prove the correctness of the algorithm and the 
min-max theorem for the original problem.

There is one more thing to check, namely that the episodes selected
for the reduced problem are distinct from the episodes selected 
at $P$ and on the $e_j$  down-edges for $j=1,2,...,\ell$. The episodes
selected in internal vertices of the $e_j$ edges are no longer
in any duplication interval in the reduced problem, and therefore
they cannot be selected there. In the reduced problem, the non-internal
vertices selected for episodes are upper endpoints of some
duplication intervals, while $P$ is no longer the upper endpoint
of any duplication interval in the reduced problem.   
$\spadesuit$

\section{Relevance for bioinformatics}\label{bioinf}

 Ohno was among the first to recognize the importance of gene and
genome
duplications \cite{ohnobook}, and the resulting opportunity for evolutionary
change afforded by
genetic redundancy.  Gene duplication, and subsequent gene loss, are the
primary drivers of changes in  gene content.
Rates of duplication and loss also have been shown to vary among lineages,
and gene loss in particular is greatly elevated
after whole-genome duplication events, which have occurred many times in the
evolution of the eukaryotes
\cite{Lynch00, NSGenetics, brunet}.
(While gene content alone has been used for reconstruction of species trees
\cite{Snel99},
it is very sensitive to parallel or convergent gains and losses, and has not
seen wide application.)

Gene trees may differ from the species tree and from each other because
 of repeated gene duplication and gene loss.
 Gene loss may
eliminate the gene from a species. A species may have more than one
representation in the gene tree as a result of gene duplication.

There are two possible explanations of finding duplicate genes: early genome duplications
and subsequent substantial gene loss, and occasional duplication of  small groups of consecutive genes, not requiring
the assumption of substantial gene loss.  The latter event is called duplication episode.
Clearly both mechanisms are present, an early vertebrate 
tetraploidization seems generally accepted.  (Duplication of "medium length" segments seem unlikely.)
 
To reconstruct the likely history of the gene content, we should know the cost associated 
with genome duplication, gene loss, and duplication episodes. We do not know those costs.  
Minimization problems for duplication episodes look for a most parsimonious explanation.

The bioinformatics literature identifies duplication vertices in the gene trees. For
every duplication vertex, 
the LCA (least common ancestor) mapping designates a vertex in the species tree, which is
a lower bound in the $\leq_S$ partial order for the point in the species tree $S$, where the
duplication of the gene could have happened. There is no absolute upper bound on this
gene duplication, but the more branchings follow the gene duplication 
 in the species tree, the more gene losses
should be assumed. Therefore a common parsimony approach is to allow the shortest duplication
interval for this gene: the edge between the lower bound vertex and its parent in the species
tree. The assumption  (\ref{treelike})   follows from the way of assigning duplication intervals
to duplication points in \cite{eulenstein} and before in the literature.


In the theorems above we have assumed that all duplication
intervals are closed and the upper and lower
endpoints of the intervals are vertices of the species tree $S$. 
Let us first consider the continuous model. Without loss of generality we may assume that
the duplication intervals end at vertices of $S$ by simply subdividing edges of
$S$ if necessary, as we did not assume that $S$ was binary. As to the assumption of
the intervals being closed: the proof really only used that the intervals are closed
upwards, i.e. they include their upper endpoint in the $>_S$ ordering.

\old{
We will view two solutions as
combinatorially equivalent if the duplication episodes they assign include the same vertices
of the species tree and use the same number of vertices on the inside of each edge;
moreover, they assign the same set of events to the vertex corresponding
to the $k$-th episode of an edge
in the $>_S$ ordering. Assume that the duplication intervals
satisfy the following property: if two duplication intervals share the same
upper endpoint $v$ end use the same edge going down from $v$, then they are either both 
open or both closed. It is easy to see that if our problem satisfies this condition, then the feasible solutions 
do not change if we cut the upper part of each interval that does not include
its $>_S$ upper endpoint by a sufficiently small $\epsilon$ and convert every interval
into one that contains its $>_S$ upper endpoint. Note in particular that this condition is
satisfied if all intervals are open upward. Moreover, if
all intervals are open upward, than an almost disjoint chain packing ${\cal C}$ is actually
a disjoint chain packing, and, since for every vertex $v$ of $S$ we have
${\cal E}(v)< 2$, the value of ${\cal C}$ is the same whether we view it as a disjoint chain
packing or an almost disjoint chain packing.
}

It is easy to see that
the discrete model can be viewed as follows: using $k$ copies of
a vertex $v$ is equivalent with placing $k$ duplication episodes on
the edge leading from $v$ to its parent (not using the parent). Therefore
in effect the discrete model is equivalent with using assuming that all duplication
intervals are of the form $[x,y)$, where $y$ is an ancestor of $x$. Thus, our
theorem for the discrete model follows from the theorem on the continuous model.


Bansal and Eulenstein \cite{eulenstein}  solved the episode minimization problem for the
discrete model. 
We leave the decision for the biology literature when the discrete or
the continuous model is to be used. We provided duality results for both models.
Sur\'anyi's Theorem~\ref{sura}  can be understood as a min-max result for the so-called {\em 
episode clustering problem}:  we are given duplication intervals and we want to
cover the duplication intervals with the minimum number of points, while we do not require 
that that duplication episodes follow in strict order.


\begin{thebibliography}{99}



\bibitem{eulenstein} M. S. Bansal, O. Eulenstein  The multiple gene
duplication problem revisited, {\em Bioinformatics}, {\bf 24}(13), (2008),
132--138.


\bibitem{speedup}  M. S. Bansal, O. Eulenstein, An $\Omega(n^2/\log n)$ speed-up
of TBR heuristics for the gene-duplication problem
{\em IEEE/ACM Transactions on Computational Biology and Bioinformatics (TCBB)}  
{\bf 5} (2008),  Issue 4,    514--524.

\bibitem{brunet}
F. Brunet, H. Crollius, M. Paris, J. Aury, P. Gibert, O. Jaillon, V. Laudet, and M. Robinson-Rechavi,  Gene loss and evolutionary rates following whole-genome duplication
 in teleost fishes, {\em Mol. Biol. and Evol.}, {\bf 23}(9) (2006), 1808--1816.

\bibitem{infer}     G. Burleigh, M. S. Bansal, O. Eulenstein, T. J. Vision,
Inferring species trees using genome duplication episodes, 
{\em BCB '10 Proceedings of the First ACM International Conference on Bioinformatics and Computational Biology}, 198--203. 

\bibitem{ilp} Wen-Chieh Chang, G.J. Burleigh, D.F. Fern\'andez-Baca, and O. Eulenstein,
An ILP solution for the gene duplication problem, {\em BMC Bioinformatics}  12(Suppl 1)(2011),  S14.
 doi:  10.1186/1471-2105-12-S1-S14



\bibitem{fellows} M. Fellows, M. Hallett, U. Stege,
On the multiple gene duplication problem, in: {\em 9th International Symposium on
Algorithms and Computation} (ISAAC '98), 
Taejon, Korea,   LNCS {\bf 1533},  1998, pp. 347--356.


\bibitem{guigo} R. Guigo {\it et al.} , Reconstruction of ancient molecular
phylogeny, Mol. Phylogenet. Evol., {\bf 6} (1996), 189--203.

\bibitem{gyarfas} A. Gy\'arf\'as, J. Lehel, A Helly type problem in trees, in: {\em
Combinatorial Theory and its Applications}, eds. P. Erd\H os, A. R\'enyi, V. T. S\'os,
North-Holland, Amsterdam, 1970, 571--584.

\bibitem{hs} A. Hajnal and J. Sur\'anyi, \"Uber die Aufl\"osung von Graphen in
vollst\"andige Teilgraphen, {\it Annales Univ. Sci. Bud. E\"otv\"os}
 {\bf 1}(1958),
115--123.

\bibitem{luo}  C. Luo, M. Chen, Y. Chen, R. Yang, H. Liu, K.  Chao, Linear-time algorithms for multiple gene duplication
problems, {\em  IEEE Transactions on Comp. Biol. and Bioinf.}, (2009).

\bibitem{Lynch00}
M. Lynch and J. Conery, The evolutionary fate and consequences of duplicate genes, {\em Science}, {\bf 290}(5494), (2000), 1151--1155.

\bibitem{NSGenetics}
J. H. Nadeau and D. Sankoff,  Comparable rates of gene loss and functional divergence after genome duplications
early in vertebrate evolution, {\em Genetics}, {\bf 147}(1997), 1259--1266.


\bibitem{ohnobook}
S. Ohno, {\em Evolution by gene duplication}, Springer-Verlag, 1970.

\bibitem{vertebrate} R.D.M. Page, J.A. Cotton, Vertebrate phylogenomics:
reconciled trees and gene duplications,  {\em Pacific Symposium on Biocomputing} (2002) 536--547.

\bibitem{Snel99}
B. Snel, P. Bork, and M. Huynen, Genome phylogeny based on gene content, {\em Nature Genetics}, {\bf 21} (1999), 108--110.



\end{thebibliography}
\end{document}